\documentclass[review]{elsarticle}

\usepackage{lineno,hyperref}
\modulolinenumbers[5]
\newtheorem{theorem}{Theorem}
\newtheorem{remark}{Remark}

\usepackage{hyperref}
\usepackage{nicematrix}
\usepackage{amsmath,blkarray}
\usepackage{makecell}
\usepackage{amssymb}
\usepackage{ulem}
\usepackage{booktabs}
\usepackage{multirow}
\usepackage{comment}
\usepackage{todonotes}

\journal{Operations Research Letters}

\begin{document}

\begin{frontmatter}

\title{A Short Proof of a Convex Representation for Stationary Distributions of Markov Chains with an Application to State Space Truncation}
% \tnotetext[mytitlenote]{Fully documented templates are available in the elsarticle package on \href{http://www.ctan.org/tex-archive/macros/latex/contrib/elsarticle}{CTAN}.}

% \author{Zeyu Zheng\fnref{myfootnote}}
% \address{Radarweg 29, Amsterdam}
% \fntext[myfootnote]{Since 1880.}

\author[mymainaddress]{Zeyu Zheng}

\author[mysecondaryaddress]{Alex Infanger}

\author[mythirdaddress]{Peter W. Glynn}

% \author[mysecondaryaddress]{Global Customer Service\corref{mycorrespondingauthor}}
% \cortext[mycorrespondingauthor]{Corresponding author}

\address[mymainaddress]{Department of Industrial Engineering and Operations Research, University of California,	Berkeley, CA 94720, USA\\}
\address[mysecondaryaddress]{Institute for Computational and Mathematical Engineering, Stanford University,
Stanford, CA 94305, USA}
\address[mythirdaddress]{Department of Management Science and Engineering, Stanford University,
Stanford, CA 94305, USA}

\begin{abstract}
In an influential paper, Courtois and Semal (1984) \cite{courtois1984bounds} establish that when $G$ is an irreducible substochastic matrix for which $\sum_{n=0}^{\infty}G^n <\infty$, then the stationary distribution of any stochastic matrix $P\ge G$ can be expressed as a convex combination of the normalized rows of $(I-G)^{-1} = \sum_{n=0}^{\infty} G^n$. In this note, we give a short proof of this result that extends the theory to the countably infinite and continuous state space settings. This result plays an important role in obtaining error bounds in algorithms involving nearly decomposable Markov chains, and also in state truncations for Markov chains. We also use the representation to establish a new total variation distance error bound for truncated Markov chains.
\end{abstract}


\end{frontmatter}

% \linenumbers

\section{Introduction}

Let $G = (G(x,dy:x,y\in S))$ be a substochastic kernel defined on a state space $S$. We say that a stochastic (transition) kernel $P=(P(x,dy):x,y \in S)$ \textit{dominates} $G$ (and write $P\ge G$) if
\begin{align*}
    P(x,B) \ge G(x,B)
\end{align*}
for all $x\in S$ and (measurable) $B \subseteq S$. A probability $\pi$ on $S$ is a \textit{stationary distribution} of $P$ if 
\begin{align}
    \pi(B) = \int_S \pi (dx) P(x,B) \label{eq:1.1} 
\end{align}
for all (measurable) subsets $B\subseteq S$. Let $\mathcal{P}$ be the set of probabilities on $S$. 

Given $G$, let 
\begin{align*}
    \Pi_G = \{ \pi \in \mathcal{P}: \pi \text{ is a stationary distribution of some} \\
    \qquad   \text{      stochastic kernel } P \text{ for which } P\ge G \}.
\end{align*}
Note that $\Pi_G$ is a convex set of probabilities. This paper provides a short proof of the following result characterizing the extreme points of $\Pi_G$. 
\begin{theorem}\label{thm::cts-state-space-Courtois-and-Semal}
Suppose that $G$ is a substochastic kernel on $S$ for which
\begin{align*}
    0 < \sum_{n=0}^{\infty} G^n(x,S) <\infty
\end{align*}
for $x\in S$. For $x\in S$, put $g(x) = \sum_{n=0}^{\infty} G^n(x,S)$, and let $\nu_x(\cdot)$ be the probability on $S$ defined by $\nu_x(\cdot) = \sum_{n=0}^{\infty} G^n(x,\cdot)/g(x) $. Then
\begin{align*}
    \Pi_G = \Big\{ \int_S \eta(dx) \nu_x(\cdot): \eta \in \mathcal{P} \Big \}. 
\end{align*}
\end{theorem}
\begin{remark}
This theorem  generalizes Theorem 2 of \cite{courtois1984bounds}, stated for finite state space, to countably infinite and continuous state space. Our theorem also does not require the irreducibility required in the 1984 theorem. However, our generalization focuses on dominating stochastic kernels, while the Courtois-Semal result permits arbitrary finite state non-negative dominating matrices. We note that strong assumptions must be imposed outside the finite state space setting to guarantee the existence of Perron-Frobenius eigenvalues, eigenfunctions, and eigenmeasures; see, for example, \cite{seneta2006non}. Our stochastic formulation avoids these issues.
\end{remark}

{
% \color{blue} 
\begin{remark}
In contrast to the proof of \cite{courtois1984bounds}, our argument is fully constructive, in the sense that for each element in each set, we construct an element in the other set that maps to that element (so neither set inclusion uses a non-constructive existence argument). 
\end{remark}
}

The above theorem plays a central role in the error bounds developed in \cite{courtois1984bounds} for approximating stationary distributions of nearly decomposable Markov chains. {\color{black} The result is also used in \cite{courtois1986computable} to develop error bounds when approximating the equilibrium distribution of a large state space Markov chain using a smaller subset (truncation) of the state space. The corresponding truncation theory for continuous time Markov jump processes can be found in \cite{dayar2011bounding}.} A survey of truncation-based algorithms for Markov jump process can be found in \cite{kuntz2021stationary}{\color{black}; see Section 4.4 and references therein for more discussion related to the above theorem.} 

% More recently, \cite{kuntz2019bounding,kuntz2021approximations} use this theorem to develop error bounds for algorithms related to state space truncations for Markov chains. 

Section 3 of this paper provides a new application of Theorem 1 to state space truncation that is exploited in \cite{glynn2022solving}. 

\section{Proof of the Theorem}
Suppose first that $P\ge G$, and that $\pi$ is a stationary distribution of $P$. Put $H(x,\cdot) = P(x,\cdot) - G(x,\cdot)$, and define
\begin{align}
    \pi_1(\cdot) &= \int_S \pi(dx) G(x,\cdot), \label{eq:2.1}\\
    \pi_2(\cdot) &= \int_S \pi(dx) H(x,\cdot). \label{eq:2.2}
\end{align}
Then, 
\begin{align}
    \pi_1(\cdot) + \pi_2(\cdot) & = \int_S \pi(dx) (G(x,\cdot)+H(x,\cdot)) \nonumber \\
    &= \int_S \pi(dx) P(x,\cdot) \nonumber\\
    &= \pi(\cdot), \label{eq:2.3}
\end{align}
since $\pi$ satisfies \eqref{eq:1.1}. Substituting \eqref{eq:2.3} into \eqref{eq:2.1}, we find that 
\begin{equation}
    \pi_1(\cdot) = \int_S \pi_1(dx) G(x,\cdot) + \kappa(\cdot), \label{eq:2.4}
\end{equation}
where
\begin{align*}
    \kappa(\cdot) = \int_S \pi_2(dx) G(x,\cdot). 
\end{align*}
Iterating \eqref{eq:2.4} $n$ times, we find that 
\begin{align*}
    \pi_1(\cdot) = \int_S \kappa(dx) \sum_{j=0}^{n} G^j(x,\cdot) + \int_S \pi_1(dx) G^{n+1}(x,\cdot).
\end{align*}
Then, \eqref{eq:2.4} implies that if we send $n\rightarrow \infty$, 
\begin{align*}
    \pi_1(\cdot) &= \int_S \kappa(dx) \sum_{j=0}^{\infty} G^j(x,\cdot)\\
    &= \int_S \pi_2(dx) \sum_{j=1}^{\infty} G^j(x,\cdot).
\end{align*}
It follows that
\begin{align}
    \pi(\cdot) &= \pi_1(\cdot) + \pi_2(\cdot) \nonumber \\
    &= \int_S \pi_2(dx) \sum_{j=0}^{\infty} G^j(x,\cdot) \nonumber \\
    &= \int_S \eta(dx)  \nu_x(\cdot) \label{eq:2.5}
\end{align}
where
\begin{align*}
    \eta(B) = \int_B \pi_2(dx) g(x).
\end{align*}

Relation \eqref{eq:2.5} implies that 
\begin{align*}
    1= \pi(S) = \int_S \eta(dx)\nu_x(S) = \eta(S),
\end{align*}
and hence $\eta$ is a probability. Hence, \eqref{eq:2.5} establishes that $\pi$ can be represented as a mixture of the $\nu_x$'s. 

Conversely, suppose that
\begin{align*}
    \mu(\cdot) = \int_S \gamma(dx) \nu_x(\cdot)
\end{align*}
for some probability $\gamma.$ Note that $g(x)\ge 1$ for $x\in S$, so 
\begin{align*}
    0<c \triangleq \int_S \gamma(dx)/g(x) <\infty.
\end{align*}
Hence, $\varphi$ as defined by
\begin{align*}
    \varphi(B) = c^{-1}\int_B \frac{\gamma(dy)}{g(y)}
\end{align*}
is a probability on $S$. 

Then, 
\begin{align}
    \mu(\cdot) &= \int_S \frac{\gamma(dx)}{g(x)} \sum_{j=0}^{\infty} G^j (x,\cdot) \nonumber\\
    &= c \int_S \varphi(dx) \sum_{j=0}^{\infty} G^j(x,\cdot). \label{eq:2.6}
\end{align}
For $x\in S$ and (measurable) $B\subseteq S$, let 
\begin{align*}
    \Phi(x,B) = \varphi(B)
\end{align*}
and set
\begin{align*}
    H(x,dy) = (1-G(x,S))\Phi(x,dy)
\end{align*}
and 
\begin{align*}
    P(x,dy) = G(x,dy) + H(x,dy)
\end{align*}
for $x,y\in S.$

Observe that \eqref{eq:2.6} implies that
\begin{align}
    (\mu P)(\cdot) &= c \int_S \varphi(dx) \sum_{j=0}^{\infty} \int_S G^j(x,dy) (G(y,\cdot) + H(y,\cdot)) \nonumber \\
    &= c\int_S \varphi(dx) \sum_{j=1}^{\infty} G^j(x,\cdot) + c\int_S \varphi(dx) \int_S \sum_{j=0}^{\infty} G^j(x,dy) H(y,\cdot) \nonumber\\
    &= \mu(\cdot) - c\varphi(\cdot) + c\int_S \varphi(dx) \int_S \sum_{j=0}^{\infty} G^j(x,dy) H(y,\cdot). \label{eq:2.7}
\end{align}
But 
\begin{align}
    & \sum_{j=0}^{\infty} \int_S  G^j(x,dy) H(y,\cdot)    \nonumber\\
    =& \sum_{j=0}^{\infty} \int_S  G^j(x,dy) (1-G(y,S))\varphi(\cdot)    \nonumber\\
        =& \sum_{j=0}^{\infty}(G^j(x,S) - G^{j+1}(x,S))\varphi(\cdot)      \nonumber\\
            =& \varphi(\cdot).    \label{eq:2.8}
\end{align}
It follows from \eqref{eq:2.7} and \eqref{eq:2.8} that
\begin{align*}
    (\mu P)(\cdot) = \mu(\cdot),
\end{align*}
proving that $\mu$ is the stationary distribution of the stochastic kernel $P$, which clearly dominates $G$. //

\section{An Application}
In this section, we briefly show how this theory applies to state space ``truncation" for Markov chain. Suppose that $Y=(Y_n:n\ge 0)$ is an $S^*$-valued Markov chain with one-step (stochastic) transition kernel $P^* = (P^*(x,dy):x,y \in S^*)$. Assume that $P^*$ has a stationary distribution $\pi^* = (\pi^*(dx):x\in S^*)$. 

Suppose that we ``truncate" the state space $S^*$ to $S\subseteq A \subseteq S^*$, and wish to compute the conditional distribution $\pi=(\pi(dx):x\in S)$ defined by
\begin{align*}
    \pi(dx) = \frac{\pi^*(dx)}{\pi^*(S)}
\end{align*}
for $x\in S.$ (This is the distribution $\pi^*$, conditioned on $S$.) If $\pi^*(S)>0$, Poincar\'{e}'s recurrence theorem (see \cite{poincare1890probleme,barreira2006poincare}) implies that $Y$ returns infinitely often a.s. to $S$ when $Y_0$ has distribution $\pi$. If $T_0 = 0$ and $T_{m+1}=\inf\{ n>T_m: Y_n \in S\}$, then $X=(X_m:m\ge 0)$, defined via $X_m = Y_{T_m}$ for $m\ge 0$ when $X_0$ has distribution $\pi$, is the \textit{Markov chain on $S$.} The Markov chain $X$ has stationary distribution $\pi$ and a transition kernel $P$ given by
\begin{align*}
    P(x,dy) = P(Y_{T_1}\in dy | Y_0 = x)
\end{align*}
for $x,y\in S$ (where $P(x,\cdot)$ can be defined via any probability on $S$ on the subset of $x$'s in $S$ having $\pi$-measure 0 for which $P(T_1<\infty|Y_0 =x) < 1$). 

To approximate $\pi$, we assume we have access only to the substochastic kernel $R=(R(x,dy):x,y\in A)$, where $R(x,dy) = P^*(x,dy)$ for $x,y\in A$. If $T=\inf\{ n\ge 0: X_n\in A^c\}$ is the first exit time from $A$, put $G(x,dy) = P(Y_{T_1}\in dy, T_1<T | Y_0=x)$ for $x,y\in S$. We note, for example, that if $A$ is a finite subset, then the substochastic kernel $G=(G(x,dy):x,y \in S)$ can be numerically computed. 

In particular, let $A' = A-S$ and set $P^*_A = (P^*_A(x,y):x,y\in A)$ with $P_A^*(x,y) = P^*(x,y)$ for $x,y \in A$. Write $P^*_A$ in block-partitioned form as 

\[
{P_A^*} = 
        \begin{blockarray}{ccc}
         & S  & A'   \\
        \begin{block}{r(rr)}
        S & P_{11} & P_{12}  \\
        A' & P_{21} & P_{22} \\
        \end{block}
    \end{blockarray} 
\]
Then, $G = P_{11}+ P_{12}(I-P_{22})^{-1}P_{21}$. 

Note that $P\ge G$, so $\pi \in \Pi_G$. It follows from Theorem 1 that
\begin{align*}
    \pi(\cdot) = \int_S \eta(dx) \nu_x(\cdot)
\end{align*}
for some $\eta \in \mathcal{P}$. Consequently, for a given ``reward" function $r:S^*\rightarrow \mathbb{R}_+$, we conclude that 
\begin{align}
    \inf_{x\in S} \int_S \nu_x(dy) r(y) \le \int_S\pi(dy) r(y) \le \sup_{x\in S} \int_S \nu_x(dy) r(y). \label{eq:3.1}
\end{align}

Furthermore, Theorem 1 implies that \eqref{eq:3.1} provides the best possible bounds on $\int_S \pi(dy) r(y)$, in the absence of auxiliary information beyond that contained in $G.$

\begin{remark}
We note that the bounds \eqref{eq:3.1} hold in general state space, even with no irreducibility assumptions on $P^*$ or $P.$
\end{remark}

\begin{remark}
The bounds \eqref{eq:3.1} play an important role in recently developed numerical truncation algorithms for discrete state space Markov chain; see \cite{glynn2022solving}.
\end{remark}

In view of Theorem 1, suppose that we choose $\tilde{\pi} = (\tilde{\pi}(dy):y\in S)$ as our approximation to $\pi$, where $\tilde{\pi}(dy) = \int_S \tilde{\eta}(dx)\nu_x(dy)$ and $\tilde{\eta}$ is some given probability on $S$. If
\begin{align*}
    \|\pi_1 - \pi_2\| \triangleq \sup_B |\pi_1(B) - \pi_2(B)|
\end{align*}
is the total variation norm between $\pi_1,\pi_2\in \mathcal{P}$, then for any $\eta_1,\eta_2\in \mathcal{P}$ and (measurable) $B \subseteq S$,  
\begin{align*}
     \int_S \eta_1(dx) \nu_x(B) - \int_S \eta_2(dy) \nu_y(B) 
    \le & \sup_{x\in S}\nu_x(B) + \sup_{y\in S}(-\nu_y(B))\\
    =& \sup_{x,y\in S} \nu_x(B) - \nu_y(B)\\
    =& \sup_{x,y\in S} |\nu_x(B) - \nu_y(B)|.
\end{align*}
Since the same inequality holds for $\eta_1$ and $\eta_2$ interchanged, it follows that, 
\begin{align}
    \|\int_S \eta_1(dx)\nu_x(\cdot) - \int_S \eta_2(dy)\nu_y(\cdot)\| \le \sup_{x,y\in S} \|\nu_x-\nu_y\|. \label{eq:3.2}
\end{align}

Since we can always specialize $\eta_1$ and $\eta_2$ to unit point mass distributions, it is evident that \eqref{eq:3.2} is the best possible total variation bound, in the sense that 
\begin{align*}
    \sup_{\eta_1,\eta_2\in\mathcal{P}} \|\int_S \eta_1(dx)\nu_x(\cdot) - \int_S \eta_2(dy)\nu_y(\cdot) \| = \sup_{x,y\in S}\|\nu_x-\nu_y\|. 
\end{align*}

In particular, we obtain a new bound on the total variation distance between the approximation $\tilde{\pi}$ and $\pi$, namely
\begin{align}
\label{eq:13}    \|\tilde{\pi} - \pi\| \le \sup_{x,y \in S} \|\nu_x - \nu_y\|. 
\end{align}

So, Theorem 1 leads to useful truncation bounds on expectations defined in terms of $\pi$, and on the total variation distance between $\tilde{\pi}$ and $\pi.$ These bounds can now be used to assess whether the choice of $S$ and $A$ achieve a given error tolerance. In particular, note that if the ``boundary layer" $A'$ is too small, there may be row vectors $\nu_x$ (with $x$ close to the boundary) that will be heavily influenced by the fact that the row $(G(x,y):y\in K)$ is very substochastic, so that $\nu_x$ looks quite different than the other $\nu_y$'s, leading to a large error bound in \eqref{eq:13}. Hence, \eqref{eq:13} can be used to choose a $K\subseteq A$ as large as possible, subject to minimization of the ``boundary layer" effect. 

% \section*{References}

% \bibliography{mybibfile}

\end{document}